\documentclass[11pt]{amsart}

\usepackage{amsfonts}
\usepackage{amsmath}
\usepackage{amssymb, esint}
\usepackage{amsthm,color}
\usepackage{enumerate}
\usepackage{mathtools}
\usepackage{cite}

\usepackage{comment}

\newcommand{\R}{{\mathbb R}}

\newcommand{\N}{{\mathbb N}}

\newcommand{\cA}{{\mathcal A}}
\newcommand{\cB}{{\mathcal B}}

\def\0{{\mathbf 0}}

\newcommand{\e}{\varepsilon}
\newcommand{\vp}{\varphi}

\newcommand{\ddiv}{\operatorname{div}}
\newcommand{\diam}{\operatorname{diam}}

\newcommand{\ra}{\rightarrow}

\def\mean#1{\mathchoice%
          {\mathop{\kern 0.2em\vrule width 0.6em height 0.69678ex depth -0.58065ex
                  \kern -0.8em \intop}\nolimits_{\kern -0.4em#1}}%
          {\mathop{\kern 0.1em\vrule width 0.5em height 0.69678ex depth -0.60387ex
                  \kern -0.6em \intop}\nolimits_{#1}}%
          {\mathop{\kern 0.1em\vrule width 0.5em height 0.69678ex
              depth -0.60387ex
                  \kern -0.6em \intop}\nolimits_{#1}}%
          {\mathop{\kern 0.1em\vrule width 0.5em height 0.69678ex depth -0.60387ex
                  \kern -0.6em \intop}\nolimits_{#1}}}

\theoremstyle{plain}
\newtheorem{thm}{Theorem}[section]

\newtheorem{cor}[thm]{Corollary}
\newtheorem{lem}[thm]{Lemma}
\newtheorem{prop}[thm]{Proposition}

\newtheorem*{claim*}{Claim}

\newtheorem{rem}[thm]{Remark}

\numberwithin{equation}{section}

\title[]{Regularity near the fixed boundary for transmission systems}

\author[A.\ Figalli]{Alessio Figalli}
\address{Department of Mathematics, ETH Z\"urich,  R\"amistrasse 101, 8092 Z\"urich, Switzerland }
\email{alessio.figalli@math.ethz.ch}

\author[S.\ Khademloo]{Somayeh Khademloo}
\address{Department of Mathematics, Faculty of Basic Sciences, Babol Noshirvani University of Technology, Babol, Iran}
\email{s.khademloo@nit.ac.ir}

\author[S.\ Kim]{Sunghan Kim}
\address{Department of Mathematics, Uppsala University, S-751 06 Uppsala, Sweden}
\email{sunghan.kim@math.uu.se}

\author[H.\ Shahgholian]{Henrik Shahgholian}
\address{Department of Mathematics, KTH Royal Institute of Technology, 100 44 Stockholm, Sweden}
\email{henriksh@kth.se}

\begin{document}

\maketitle

\begin{abstract}
Given $\Omega\subset\R^n$ with $n\geq 2$, $D\subset \Omega$ open, and $u:\Omega \to \R^m$, we study elliptic systems of the type
$$ \ddiv \big( ( \cA + (\cB- \cA)\chi_D)\nabla u\big) = 0 \quad \text{in $\Omega\cap B_1$,} $$
for some uniformly elliptic tensors $\cA$ and $\cB$ with
H\"{o}lder continuous entries. We show that, given appropriate boundary data, the Lipschitz regularity of $u$ inside $B_1 \cap D$ is transmitted to $B_{1/2}\cap \Omega$ up to the boundary of $\Omega$.
This corresponds to the boundary counterpart of the interior regularity results in \cite{FKS22}.

\bigskip 

\begin{center}
{\it In honor of Nina Nikolaevna Uraltseva for her 90th birthday.}
\end{center}


\end{abstract}

\tableofcontents


\section{Introduction}
\subsection{Background and motivation}
Transmission problems describe phenomena in which a physical quantity changes behavior across some fixed surface that we call \textit{interface}. In mathematical physics, these problems involve interfaces
immersed in material bodies that contain two or more components with distinct physical characters. Any physical processes across
the interfaces could be \textit{interrupted}
and lose some continuity. So, transmission problems deal with a fixed interface where solutions change \textit{abruptly}.
 The study of their behavior across this interface is of great interest these days. Mathematically, the transmission problems are described by PDEs on each individual component, and then the solutions are glued together through the so-called transmission conditions imposed on the interfaces. Such problems cannot be treated as the usual boundary value problems because the solutions in touching components will
interact with each other by the transmission conditions, and so there is noticeable contrast between them in essence and methods.  The analysis of such problems
started with the pioneering work of Picone \cite{Pic54}. Picone, for the first time in 1954, introduced transmission problems in classical elasticity theory.
 The transmission problem for elliptic equations with smooth coefficients and interfaces was developed by Schechter in \cite{Sch60}.
 In fact, he generalized the theory to include smooth elliptic operators in non-divergence form in domains with smooth interfaces, and since
then transmission problems have been of great interest due to their applications in different areas such as electromagnetic processes, composite materials,
vibrating folded membranes, climatology, etc.; see, for instance, \cite{Bor10}.
The main issues for transmission problems are optimal interior regularity of solutions and regularity of solutions near the interfaces, which is expected to be lower than for a usual PDE without transmission, due to its intrinsic nature.

In this paper, we study the boundary regularity of a transmission problem where Lipschitz continuity is assumed in an {\it a priori} unknown inclusion $D$. The interior case is studied in \cite{FKS22}, where the authors showed the propagation of Lipschitz regularity in a uniform neighborhood of the inclusion. Here, we aim at extending the uniform propagation result all the way up to the fixed boundary layer.

\subsection{Standing assumptions and main results}

Let $\Omega$ be a domain in $\R^n$ with $0\in\partial\Omega$, $D\subset \Omega$ open, and $B_1$ the unit  ball centered at the origin. We consider the following transmission problem
\begin{equation}\label{eq:main}
 \begin{cases}   
 \ddiv \big( ( \cA + (\cB- \cA)\chi_D)\nabla u\big) = 0 & \text{in }B_1 \cap\Omega,\\
u = g &\text{on }B_1 \cap\partial\Omega,
\end{cases}
\end{equation}
where $\cA = (A_{ij}^{\alpha\beta})_{1\leq i,j\leq m}^{1\leq\alpha,\beta\leq n}$ and $\cB= (B_{ij}^{\alpha\beta})_{1\leq i,j\leq m}^{1\leq\alpha,\beta\leq n}$ are the coefficient tensors, and $g = (g^1,\cdots,g^m)$ is a boundary datum. We say $u\in W^{1,2}(B_1 \cap\Omega;\R^m)$ is a weak solution of \eqref{eq:main}, if 
$$
\int_{B_1 \cap\Omega} \Big(A_{ij}^{\alpha\beta} + \Big(B_{ij}^{\alpha\beta} - A_{ij}^{\alpha\beta}\Big)\chi_D\Big) \partial_\beta u^j\partial_\alpha \vp^i \,dx = 0,
$$
for every $\vp\in W_0^{1,2}(B_1 \cap\Omega;\R^m)$ and verifies $u|_{B_1 \cap\partial\Omega} =g$ in the trace sense. Here and thereafter, we shall use the summation convention for repeated indices. 

Throughout this paper, $n\geq 2$ is the dimension of the ambient space, $m\geq 1$ that of the target space, and we fix constants $\lambda\in(0,1)$, 
$\sigma\in(0,1)$, and $K > 0$. Let us specify the assumptions on $\cA$, $\cB$ and $g$ as follows:\\
$\bullet$ (Strong ellipticity) For every $\xi\in\R^{mn}$,
    \begin{equation}\label{eq:ellip}
    \lambda |\xi|^2 \leq A_{ij}^{\alpha\beta}\xi_\alpha^i\xi_\beta^j, B_{ij}^{\alpha\beta}\xi_\alpha^i\xi_\beta^j\leq \frac{1}{\lambda}|\xi|^2\qquad \text{  a.e.\ in $B_1 \cap\Omega$.}
    \end{equation}
$\bullet$ (Regularity) For every $i,j,\alpha,\beta$ it holds $A_{ij}^{\alpha\beta} \in C^{0,\sigma}(B_1 \cap\Omega)$, with the estimate
    \begin{equation}\label{eq:A-Ca}
    [A_{ij}^{\alpha\beta}]_{C^{0,\sigma}(B_1 \cap\Omega)}\leq K.
    \end{equation}

Our main result reads as follows.

\begin{thm}\label{thm:main}
Let $D\subset\Omega$ be an open set.
Assume that $B_1\cap\partial\Omega$ is of class $C^{1,\sigma}$, and $g\in C^{1,\sigma}(B_1 \cap\partial\Omega;\R^m)$. Let $\cA,\cB: B_1 \cap\Omega\to\R^{n^2m^2}$ satisfy \eqref{eq:ellip} and \eqref{eq:A-Ca}, and $u\in W^{1,2}(B_1 \cap\Omega;\R^m)$ be a weak solution to \eqref{eq:main}. Assume that  $|\nabla u| \in L^\infty(B_1 \cap D)$. Then we have $|\nabla u| \in L^\infty(B_{1/2} \cap \Omega)$ with the estimate  
$$
\| \nabla u \|_{L^\infty(B_{1/2} \cap\Omega)} \leq C\bigg( \bigg(\mean{B_1 \cap\Omega}|\nabla u|^2\,dx\bigg)^{1/2} + \| \nabla u \|_{L^\infty(B_1 \cap D)}  + \| \partial_\tau g \|_{C^{0,\sigma}(B_1 \cap\partial\Omega)} \bigg), 
$$
where $C$ depends only on $n$, $m$, $\lambda$, $K$, the $C^{1,\sigma}$-character of $B_1 \cap\partial\Omega$; here $\partial_\tau$ is the tangential derivatives on $B_1 \cap\partial\Omega$.
\end{thm}

It is noteworthy that our proof works similarly well when  $0 \in  \partial D\cap \partial\Omega \neq\emptyset$, or generally if   $B_1 \cap \partial D\cap \partial\Omega \neq\emptyset$.

Combining Theorem~\ref{thm:main} with \cite[Theorem 1.1]{FKS22}, we obtain the global Lipschitz estimate for weak solutions to elliptic transmission problems.

\begin{cor}\label{cor:main}
    Let $\Omega$ be a bounded $C^{1,\sigma}$-domain, $D\subset\Omega$ be an open set, $g\in C^{1,\sigma}(\overline\Omega;\R^m)$, and $\cA,\cB:\Omega\to\R^{n^2m^2}$ satisfy \eqref{eq:ellip} and \eqref{eq:A-Ca} in the whole domain $\Omega$. Let $u\in g + W_0^{1,2}(\Omega;\R^m)$ be a weak solution to \eqref{eq:main}. If $|\nabla u|\in L^\infty(D)$, then $|\nabla u|\in L^\infty(\Omega)$, and 
    $$
    \| \nabla u \|_{L^\infty(\Omega)} \leq C \bigg( \| g\|_{C^{1,\sigma}(\partial\Omega)} + \| \nabla u \|_{L^\infty(D)} \bigg),
    $$
    where $C$ depends only on $n$, $m$, $\lambda$, $K$, the $C^{1,\sigma}$-character of $\partial\Omega$, and $\diam\Omega$. 
\end{cor}


\section{Reduction to the flat  boundary case}\label{sec:reduce}

Throughout this paper we  set $B_1^+ := B_1\cap\{x: x_n >0\}$ and $B_1' := B_1\cap\{x: x_n = 0\}$. 

As we assume that $B_1 \cap\partial\Omega$ is of class $C^{1,\sigma}$, there is a boundary flattening map $\Phi:B_1 \cap\overline\Omega \to B_1^+\cup B_1'$ with $\Phi\in C^{1,\sigma}(B\cap\overline\Omega)$. Through the transformation
$$
A_{ij}^{\alpha\beta} \mapsto 
\frac{D_\gamma\Phi^{\alpha}\circ\Phi^{-1} D_\delta\Phi^{\beta}\circ \Phi^{-1} A_{ij}^{\gamma\delta}\circ\Phi^{-1}}{|\det (D\Phi\circ\Phi^{-1})|}
$$
and analogously for $B_{ij}^{\alpha\beta}$, 
it suffices to deal with the system \eqref{eq:main} in $B_1^+\cup B_1'$ in place of $B_1 \cap\overline\Omega$. 

Note that the above transformation does not alter the structural conditions \eqref{eq:ellip} and \eqref{eq:A-Ca}, modulo slight modifications in the involved parameters $\lambda$ and $K$; the changes depend only on $n$, $m$, the $C^{1,\sigma}$-character of $B_1 \cap\partial\Omega$.

For this reason, from now on, we shall assume that $B_1 \cap\Omega = B_1^+$ and $B_1 \cap\partial\Omega = B_1'$. We denote by $\nabla'$ the gradient in the first $(n-1)$ variables. Rescaling the domain and the target, our main theorem follows from the following proposition via a standard covering argument. 

\begin{prop}\label{prop:main}
    Assume that $\cA, \cB: B_1^+\to\R^{n^2m^2}$ verify \eqref{eq:ellip} and $g \in C^{1,\sigma}(B_1';\R^m)$. Let $u\in W^{1,2}(B_1^+;\R^m)$ be a weak solution to \eqref{eq:main}. Then there exist $\e_0 > 0$ and $c>1$, depending only on $n$, $m$, $\lambda$ and $\sigma$, such that if 
\begin{equation}\label{eq:u}
\|\nabla u\|_{L^2(B_1^+)} \leq 1,\qquad \|\nabla u \|_{L^\infty(B_1^+\cap D)} \leq \e_0,
\end{equation}
\begin{equation}\label{eq:A-g}
\| \nabla' g \|_{C^{0,\sigma}(B_1')} \leq \e_0,\qquad [A_{ij}^{\alpha\beta}]_{C^{0,\sigma}(B_1^+)} \leq \e_0, 
\end{equation}
then $|\nabla u| \in L^\infty(B_{1/2}^+)$ with 
\begin{equation}\label{eq:C01-re}
\| \nabla u \|_{L^\infty(B_{1/2}^+)} \leq c.
\end{equation}
\end{prop}

In the rest of the paper, we shall call a constant universal if it depends only on the fixed parameters $n$, $m$, $\lambda$, and $\sigma$. Also, by $c$ and $C$ we shall denote positive universal constants, which may differ at each occurrence. Unless stated otherwise, we shall also assume that $\cA$, $\cB$, $g$ verify \eqref{eq:ellip} and \eqref{eq:A-g}, and $u$ is a weak solution to \eqref{eq:main} satisfying \eqref{eq:u}. 


\section{$BMO$-type estimate}\label{sec:bmo}

 This section is devoted to the proof of a universal $BMO$-bound of $\nabla u$ up to the boundary. This is mainly due to the assumption that $|\nabla u|$ is bounded inside the unknown inclusion $D$. The key here is that the approximating affine mappings at each scale can be chosen to be constant in the first $(n-1)$ variables. 
\begin{lem}\label{lem:BMO}
There exist $\e_0 > 0$ and $C>1$ universal such that if \eqref{eq:u} and \eqref{eq:A-g} hold, then for every $z=(z',z_n)\in B_{1/2}^+\cup B_{1/2}'$ and every $r\in (0,1/4)$ one can find an affine map $\ell_{z,r}$ such that 
\begin{equation}\label{eq:Du-BMO}
\int_{B_r(z)\cap B_1^+} |\nabla (u - \ell_{z,r})|^2\,dx \leq Cr^n. 
\end{equation}
Also, for $r > z_n$ one can choose $\ell_{z,r}=\ell_{z',2r}$ such that $\nabla'\ell_{z,r} = 0$. 
\end{lem}

\begin{proof}
By \eqref{eq:A-g} we can extend $g$ to $B_1^+$ so that $g\in C^{1,\sigma}(B_1^+)$ and 
\begin{equation}\label{eq:g-ext-C1a}
\| \nabla g \|_{C^{0,\sigma}(B_1^+)} \leq c \| \nabla' g \|_{C^{0,\sigma}(B_1')} \leq c\e_0. 
\end{equation}
Consider the auxiliary map $v: = u - g$, which solves 
$$
\begin{cases}
    \ddiv ( \cA \nabla v) = \ddiv F & \text{in }B_1^+ \\
    v = 0 &\text{on }B_1'
\end{cases}
$$
in the weak sense, where $F := (\cA- \cB)\chi_D \nabla u  - \cA\nabla g$. Due to \eqref{eq:ellip}, \eqref{eq:u}, \eqref{eq:A-g}, and \eqref{eq:g-ext-C1a}, we have 
\begin{equation}\label{eq:F-Linf}
\| F \|_{L^\infty(B_1^+)}\leq c\e_0.
\end{equation}
Choosing $\e_0$ smaller if necessary, \eqref{eq:u} and \eqref{eq:g-ext-C1a} imply that
\begin{equation}\label{eq:Dv-L2}
\| \nabla v \|_{L^2(B_1^+)} \leq 1 + c\e_0 \leq 2. 
\end{equation}

Hence, a standard\footnote{This boils down to a compactness argument followed up with an iteration. See Appendix in \cite{FS15} for a similar line of arguments for obstacle problem.   }
approximation argument shows the existence of a universal constant $\e_0>0$ such that if \eqref{eq:A-g} and \eqref{eq:F-Linf} hold with $\e_0$, then for every $z' \in B_{1/2}'$ and every $r\in(0,1-|z'|)$ there exists an affine map $\ell_{z',r}$ such that $\nabla' \ell_{z',r} = 0$ and 
$$
\int_{B_r^+(z')} |\nabla (v - \ell_{z',r})|^2\,dx \leq Cr^n. 
$$
Then for any $z\in B_{1/2}^+\cup B_{1/2}'$ with $z = (z',z_n)$, we choose $\ell_{z,r} := \ell_{z',2r}$ for every $r\in (z_n,\frac{1}{4})$. Since in this case $B_r(z) \cap B_1^+\subset B_{2r}(z')$, we deduce that 
\begin{equation}\label{eq:Dv-BMO}
\int_{B_r(z)\cap B_1^+} |\nabla (v - \ell_{z,r})|^2\,dx \leq C(2r)^n. 
\end{equation}
Recalling that $\nabla' \ell_{z,r} = \nabla' \ell_{z',2r} = 0$,
thanks to \eqref{eq:g-ext-C1a} and \eqref{eq:Dv-BMO} we get  
$$
\int_{B_r(z)\cap B_1^+} |\nabla(u-\ell_{z,r})|^2\,dx \leq 2 c(2r)^n + 2\int_{B_r(z)\cap B_1^+} |\nabla g|^2\,dx \leq Cr^n,
$$
as desired. 

Finally, for $r\in(0,z_n)$, we can follow \cite[Proof of Step 1B]{Acq92} to obtain \eqref{eq:Dv-BMO}. Then, in the same way as above, we deduce \eqref{eq:Du-BMO} via the triangle inequality and \eqref{eq:g-ext-C1a}. 
\end{proof}




As a corollary, we obtain a bound on the oscillation of $\nabla\ell_{z,r}$ as a function of $r$.
\begin{cor}\label{cor:BMO}
Let $\ell_{z,r}$ be as in Lemma \ref{lem:BMO}. Then there exists a universal constant $C$ such that for all $r\in (z_n,\frac{1}{4})$,
\begin{equation}\label{eq:lr}
|\nabla \ell_{z,r} -\nabla \ell_{z,r/2}|\leq C.
\end{equation}
\end{cor}

\begin{proof}
Using estimate \eqref{eq:Du-BMO}, one has
\begin{equation*}
\mean{B_r(z)\cap B_1^+}|\nabla u-\nabla \ell_{z,r}|^2\,dx+ \mean{B_{r/2}(z)\cap B_1^+} |\nabla u-\nabla \ell_{z,r/2}|^2\,dx\leq C.
\end{equation*}
Therefore, we obtain 
$$
\begin{aligned}
  &|\nabla \ell_{z,r}-\nabla \ell_{z,r/2}|^2 \\
  &=\mean{B_r(z)\cap B_1^+}|\nabla \ell_{z, r}-\nabla \ell_{z, r/2}|^2dx \\
   &\leq 2\mean{ B_r (z) \cap B_1^+ } |\nabla u-\nabla \ell_{z, r}|^2\,dx+2\mean{B_r (z) \cap B_1^+ } |\nabla u-\nabla \ell_{z,r/2}|^2\,dx\\
   &\leq  2\mean{B_r (z) \cap B_1^+  } |\nabla u-\nabla \ell_{z, r}|^2\,dx+2^{1-n}\mean{ B_{r/2} (z)  \cap B_1^+ } |\nabla u-\nabla \ell_{z, r/2}|^2\,dx\\
   &\leq C,
\end{aligned}
$$
for some universal constant $C$.
\end{proof}


\section{Decay of the density}\label{sec:decay}

Throughout this section, we shall assume that $\cA$, $\cB$, $g$ and $u$ verify \eqref{eq:ellip}, \eqref{eq:A-g}, and \eqref{eq:u}, with $\e_0$ chosen as in Lemma \ref{lem:BMO}. We write
$$
D_{z,r} := \{ y : z + ry \in B_1^+\cap D\},\qquad \Lambda_{z,r} := \frac{| B_r(z)\cap B_1^+\cap D|}{|B_r(z)\cap B_1^+|}.
$$
We follow the same idea as in the interior case, see \cite[Lemma 2.1]{FKS22}. Note that the size of $|\nabla u|$ at scale $r$ is comparable with $|\nabla \ell_r|$
 by Lemma \ref{lem:BMO}. Hence, the lemma below shows that if the size of $|\nabla u|$ at scale $r$ is large enough, then the density  $\Lambda_ {z,r}$
 must have a universal decay from  $r$ to $r/2$.

\begin{lem}\label{lem:decay}
Let $z \in B_{1/2}^+$ and $r\in(z_n,1/4)$ be given, and let $\ell_{z,r}$ be as in Lemma \ref{lem:BMO}. Given any $\mu \in (0,1)$, there exist constants $M(\mu) , C(\mu )> 1$, both depending only on $\mu$ and the universal parameters, such that if $|\nabla\ell_{z,r}| \geq M(\mu)$ then
$$
\Lambda_{z',r/2} \leq \mu\Lambda_{z',r} + C(\mu) r^{3n\sigma}.
$$
\end{lem}

\begin{proof}
Throughout this proof, $C > 1$ will denote a universal constant that may differ at each occurrence. Without loss of generality we  shall present the proof only for the case $z' = 0$.

To simplify the notation, we denote by $\ell_r$, $D_r$, $\Lambda_r$ the quantities $\ell_{0,r}$, $D_{0,r}$, $\Lambda_{0,r}$. It follows from Lemma \ref{lem:BMO} that, since $z' =0$, for $r>z_n$ we can choose $\ell_{z,r}$ to be equal to $\ell_{2r}$. In view of Corollary \ref{cor:BMO}, we may also assume that $|\nabla \ell_r| \geq \frac{1}{2}M(\mu)$ (by, say, taking $M(\mu) > 2C$ with $C$ as in \eqref{eq:lr}).

Set $u_r(y)=u(ry)/r$, $g_r(y) := g(ry)/r$, and onsider an auxiliary problem, 
\begin{equation}\label{eq:v}
\begin{cases}
\ddiv(\cA(0)\nabla v_r)=0 \quad\text{in }B_1^+,\\
v_r - (u_r - \ell_r) \in  W_0^{1,2}(B_1^+),
\end{cases}
\end{equation}
which admits a unique weak solution, due to \eqref{eq:ellip}. By \eqref{eq:ellip} again and \eqref{eq:Du-BMO}, we also have the energy estimate, 
\begin{equation}\label{eq:Dvr-L2}
\| \nabla v_r \|_{L^2(B_1^+)}\leq C||\nabla(u_r-\ell_r)||_{L^2 (B_{1}^+)}  \leq C.
\end{equation}
Since $r > z_n$, Lemma \ref{lem:BMO} ensures that $|\nabla'\ell_r| = 0$. Hence, by \eqref{eq:A-g}, 
\begin{equation}\label{eq:gr-lr-C1a}
\| \nabla' (g_r - \ell_r) \|_{C^{0,\sigma}(B_1')} \leq \e_0. 
\end{equation}
By \eqref{eq:ellip}, \eqref{eq:Dvr-L2}, and \eqref{eq:gr-lr-C1a},
applying the standard up-to-boundary gradient estimate to \eqref{eq:v} yields 
\begin{equation}\label{es:v}
\| \nabla v_r \|_{L^\infty(B_{3/4}^+)} \leq C.
\end{equation}
As a byproduct of \eqref{es:v}, we have
\begin{equation}\label{es:v_r on D_r}
\|\nabla v_r\|_{L^p(B_{3/4}^+\cap D_r)}\leq C|D_r|^{1/p},
\end{equation}
for every $ p > 1$.

The rest of the proof follows the same lines of \cite[Lemma 2.1]{FKS22}. We only need to apply the boundary $W^{1,p}$ estimates instead of the interior $W^{1,p}$ estimates used there. Here, we present the argument only to keep our writing self-contained. 

Set $w_r := u_r - \ell_r - v_r \in W_0^{1,2}(B_1^+;\R^m)$, which solves 
\begin{equation}\label{eq:w}
\ddiv(\cA_r\nabla w_r)=\ddiv (F_r + G_r) \quad \text{in }\Omega_r
\end{equation}
in the weak sense, where we set 
\begin{equation}\label{def:F_r}
\begin{split}
&\cA_r (x) := \cA(rx), \qquad G_r(x) := (\cA_r(0) - \cA_r(x))(\nabla\ell_r + \nabla v_r) \\
&\cB_r (x) := \cB(rx), \qquad F_r(x):=(\cA_r(x)- \cB_r(x))\chi_{D_r}\nabla u_r. 
\end{split}
\end{equation}
By  \eqref{eq:main} and \eqref{eq:ellip}, for any $p > 1$ we have
\begin{equation}\label{F_r}
\begin{aligned}
\int_{B_1^+}| F_r|^p\,dx &\leq C^p |D_r| .
\end{aligned}
\end{equation}
Next, by \eqref{eq:A-g} and \eqref{es:v} (as well as our assumption that $|\nabla \ell_r| \geq M (\mu) > 1/2$), 
 \begin{equation}\label{G_r}
 ||G_r||_{L^\infty(B_{3/4}^+)}^2 + \int_{B_{3/4}^+} |G_r|^2\,dx  \leq Cr^{2\sigma}|\nabla \ell_r|^2.
 \end{equation}
By \eqref{eq:ellip}, \eqref{eq:A-g},  \eqref{F_r}, and \eqref{G_r}, the up-to-the-boundary $W^{1,p}$-estimate applied to \eqref{eq:w} yields that, for every $p>1$, 
\begin{equation}\label{w_r}
\int_{B_{1/2}^+}|\nabla w_r|^p\,dx \leq C_p\big(|D_r| + r^{p\sigma}|\nabla \ell_r|^p\big).
\end{equation}
Finally, we proceed as follows: since $\ell_r \equiv u_r - v_r - w_r$ in $\Omega_r$, and $\nabla\ell_r$ is constant (as $\ell_r$ is linear), we obtain 
$$
\begin{aligned}
|B_{1/2}^+\cap D_r||\nabla\ell_r|^p &= \int_{B_{1/2}^+\cap D_r} |\nabla(u_r - v_r - w_r)|^p\,dx \\
& \leq 3^p \int_{B_{1/2}^+\cap D_r} (|\nabla u_r|^p + |\nabla v_r|^p + |\nabla w_r|^p) \,dx \\
& \leq C_p \big(|D_r| + r^{p\sigma} |\nabla\ell_r|^p\big). 
\end{aligned}
$$
Note that $|B_{1/2}^+\cap D_r| = 2^{-n}|D_{r/2}|$. Therefore, choosing $p = 3n$, we can choose a constant $M(\mu)> 1$, depending on the given constant $\mu$ and the universal parameters, such that $|\nabla\ell_r| \geq \frac{1}{2}M(\mu)$ implies 
$$
\Lambda_{r/2} \leq \mu\Lambda_r + C(\mu) r^{3n\sigma},
$$
as desired. 
\end{proof}

\section{Proof of Theorem \ref{thm:main}}\label{sec:pf-main}

As noted in Section \ref{sec:reduce}, Theorem \ref{thm:main} follows directly via boundary flattening argument from Proposition \ref{prop:main}. Hence, we shall focus on the proof of the latter. 

\begin{proof}[Proof of Proposition \ref{prop:main}]

We follow the dichotomy argument in \cite[Theorem 1.1]{FKS22} with suitable adjustments due to the presence of the boundary; it should be stressed that the idea originated from \cite{ALS13} and \cite{FS14}. Nevertheless, we need an extra step to take care of the boundary here. 

It suffices to consider $z \in B_{1/2}^+$ a Lebesgue point $z$ of $|\nabla u|^2$; since $|\nabla u|^2\in L^1(B_1^+)$, almost every point in the cylinder is a Lebesgue point. We split the proof into two cases: 
\begin{enumerate}[({Case} 1)]
\item $\liminf\limits_{k\ra\infty}  \left| \nabla \ell_{z,2^{-k}} \right| < 3M$,
\item $\liminf\limits_{k\ra\infty}  \left| \nabla \ell_{z,2^{-k}} \right| \geq 3M$,
\end{enumerate}
where $M > 1$ is the large, universal constant from Lemma \ref{lem:decay}. More precisely, we choose $\mu \in (0,1)$ to be a small universal constant such that
\begin{equation}\label{eq:mu}
2^{3n\sigma}\mu \leq 1,
\end{equation}
and take $M \equiv M(\mu)$ as in the lemma. To further simplify our notation, assume (without losing any generality) that $z = (0,z_n)$.

As for Case 1, \eqref{eq:Du-BMO} yields 
$$
\liminf_{r\to 0^+} \mean{B_r(z)\cap B_1^+} |\nabla u|^2 \,dx \leq C. 
$$
However, $z$ being a Lebesgue point of $|\nabla u|^2$, the bound for the limit infimum of the average integral is enough to deduce that $|\nabla u(z)|\leq C$. 

As for Case 2, we define
\begin{equation}\label{k_0}
k_0=\inf\{k\geq 0: |\nabla\ell_{z,2^{-j}}|\geq M\quad \forall j\geq k\}.
\end{equation}
Following the lines of the proof of \cite[Theorem 1.1]{FKS22}, by taking $M > 1$ large, we see that $k_0$ is a finite integer with $k_0\geq 2$, and that 
\begin{equation}\label{eq:k0}
|\nabla\ell_{z,2^{-k_0}}|\leq 2M. 
\end{equation}
If $2^{-k_0} \leq z_n$, then $B_{2^{-k_0}}(z)\subset B_1^+$. Therefore, \eqref{eq:Du-BMO} and \eqref{eq:k0} imply 
\begin{equation}\label{eq:Du-L2-k0}
\mean{B_{2^{-k_0}}(z)} |\nabla u|^2\,dx \leq C.
\end{equation}
As $B_{2^{-k_0}}(z)\subset B_1^+$, we can employ the interior estimate from \cite[Theorem 1.1]{FKS22} to $u$. So, by \eqref{eq:Du-L2-k0} and \eqref{eq:u}, as well as the choice of $z$ as a Lebesgue point,  we conclude that
$$
|\nabla u(z)| \leq C.  
$$
Thus, we are only left with the case $2^{-k_0} > z_n$. Choose $k_1$ such that 
\begin{equation}\label{eq:k1}
k_1 := \max\{ k \geq k_0: 2^{-k} > z_n\}.
\end{equation}
By \eqref{eq:mu}, \eqref{eq:k0}, and \eqref{eq:k1}, we can iterate Lemma \ref{lem:decay} to obtain (recall that $z'=0$) 
\begin{equation}\label{eq:Lk}
\Lambda_{2^{-{k_0 + j}}} \leq C 2^{-nj\sigma}\qquad\forall\, j\in [0,k_1 - k_0]\cap\N. 
\end{equation}
Let $u_r$, $g_r$, $\cA_r$, and $F_r$ be as in the proof of Lemma \ref{lem:decay}, and set $\hat G_r := (\cA(0) - \cA_r)\nabla\ell_r$. Define $\hat w_{2^{-k_0}} := u_{2^{-k_0}}-\ell_{2^{-k_0}}$ (note that $\ell_{2^{-k_0}} = \ell_{0,2^{-k_0}}$ which may not necessarily be equal to $\ell_{z,2^{-k_0}}$), which solves
\begin{equation}\label{eq:hat{w}}
\begin{cases}
\ddiv(\cA_{2^{-k_0}}\nabla \hat w_{2^{-k_0}})=\ddiv (F_{2^{-k_0}} + \hat G_{2^{-k_0}}) & \text{in }B_1^+\\
\hat w_{2^{-k_0}} = g_{2^{-k_0}}-\ell_{2^{-k_0}} & \text{on }B_1'
\end{cases}
\end{equation}
in the weak sense. By \eqref{eq:Lk} and \eqref{eq:u}, for each $j=0,1,\cdots,k_0-k_1$ we have
\begin{equation}\label{eq:Fk-L2}
\mean{B_{2^{-j}}^+} |F_{2^{-k_0}}|^2\,dy \leq C\e_0^2 2^{-nj\sigma}.
\end{equation}
On the other hand, by \eqref{eq:A-g}, \eqref{eq:k0}, and \eqref{eq:lr}, we have 
\begin{equation}\label{eq:Gk-L2}
\mean{B_{2^{-j}}^+} |\hat G_{2^{-k_0}}|^2\,dy \leq CM^2\e_0^2 2^{-(k_0 + j)n\sigma}.
\end{equation}
Moreover, by \eqref{eq:Du-BMO}, 
\begin{equation}\label{eq:Dhwk0-L2}
\int_{B_1^+} |\nabla\hat w_{2^{-k_0}}|^2\,dy \leq C. 
\end{equation}
Furthermore, since Lemma \ref{lem:BMO} ensures $|\nabla'\ell_{2^{-k_0}} | = 0$, we have from \eqref{eq:hat{w}} that $\nabla' w_{2^{-k_0}}|_{B_1'} = \nabla' g_{2^{-k_0}}$. Thus, by \eqref{eq:A-g}, \eqref{eq:Fk-L2}, \eqref{eq:Gk-L2}, and \eqref{eq:Dhwk0-L2}, we can employ a standard approximation argument to deduce that 
\begin{equation}\label{eq:Dhwk0}
    \mean{B_{2^{-j}}^+} |\nabla \hat w_{2^{-k_0}}|^2\,dy \leq C 2^{-nj\sigma} \mean{B_1^+} |\nabla \hat w_{2^{-k_0}}|^2\,dy \leq C2^{-nj\sigma},
\end{equation}
for every integer $k =0,1,\cdots,k_1-k_0$, where the last inequality follows from \eqref{eq:Du-BMO} after rescaling. In particular, taking $j = k_1 - k_0$ in \eqref{eq:Dhwk0}, and utilizing the triangle inequality as well as \eqref{eq:k0}, we deduce that 
\begin{equation}\label{eq:Duk0}
\begin{aligned}
\mean{B_{2^{-k_1}}^+}|\nabla u|^2\,dx & = \mean{B_{2^{-(k_1 - k_0)}}^+} |\nabla u_{2^{-k_0}}|^2\,dy \\
&\leq 2 \mean{B_{2^{-(k_1- k_0)}}^+} (|\nabla \hat w_{2^{-k_0}}|^2 + |\nabla \ell_{2^{-k_0}}|^2) \,dy \leq C. 
\end{aligned}
\end{equation}
Finally, recall from \eqref{eq:k1} that $2^{-k_1-1}\leq z_n < 2^{-k_1}$. This implies that $B_{2^{-k_1-1}}(z)\subset B_{z_n}(z) \subset B_{2^{-k_1}}^+$. Whence, by \eqref{eq:Duk0},
$$
\mean{B_{z_n}(z)} |\nabla u|^2\,dx \leq C\frac{|B_{2^{-k_1}}^+|}{|B_{z_n}(z)|} \leq C. 
$$
Now we can argue as at the beginning of Case 2: applying \cite[Theorem 1.1]{FKS22} to $u$ in $B_{z_n}(z)$, with the help of \eqref{eq:Duk0} and \eqref{eq:u}, as well as the choice of $z$ as a Lebesgue point, we get
$$
|\nabla u(z)| \leq C.  
$$
This covers all possible cases, completing the proof.
\end{proof}

\begin{rem}[Parabolic counterpart]
Following the approach outlined in \cite{FKS22} and utilizing the methodology described in the preceding section, one could explore the regularity of linear transmission parabolic systems. This line of inquiry is not pursued in the current paper, and we defer this extension to future research endeavors. 
\end{rem}




 \section*{Acknowledgements}
Sunghan Kim was supported by a grant from the Verge Foundation. Henrik Shahgholian was supported by the Swedish Research Council (grant no.~2021-03700).

\section*{Declarations}

\noindent {\bf  Data availability statement:} All data needed are contained in the manuscript.

\medskip
\noindent {\bf  Funding and/or Conflicts of interests/Competing interests:} The authors declare that there are no financial, competing or conflict of interests.

\end{document}